\documentclass[a4paper]{amsart}
\usepackage[all]{xy}
\usepackage{ifthen}
\usepackage{amssymb}

\SelectTips{eu}{}
\calclayout

\title{Expansions of abelian categories}

\thanks{Chen was supported by the
  Alexander von Humboldt Stiftung and the National Natural Science
  Foundation of China (No.\ 10971206).}

\author[Xiao-Wu Chen]{Xiao-Wu Chen}
\address{Xiao-Wu Chen\\ Department of Mathematics\\ University of Science and
Technology of China\\ Hefei 230026, Anhui\\ PR China.}
\email{xwchen@ustc.edu.cn}

\author[Henning Krause]{Henning Krause}
\address{Henning Krause\\ Fakult\"at f\"ur Mathematik\\
Universit\"at Bielefeld\\ D-33501 Bielefeld\\ Germany.}
\email{hkrause@math.uni-bielefeld.de}

\numberwithin{equation}{subsection}

\theoremstyle{plain}
\newtheorem{lem}[equation]{Lemma}
\newtheorem{prop}[equation]{Proposition}

\newtheorem{thm}[equation]{Theorem}
\newtheorem*{Thm}{Theorem}

\theoremstyle{remark}

\theoremstyle{definition}
\newtheorem{rem}[equation]{Remark}
\newtheorem{exm}[equation]{Example}
\newtheorem{defn}[equation]{Definition}

\DeclareMathOperator{\add}{add}

\DeclareMathOperator{\coh}{coh}

\DeclareMathOperator{\End}{End}

\DeclareMathOperator{\Hom}{Hom}
\DeclareMathOperator{\Ext}{Ext}

\DeclareMathOperator{\Ker}{Ker}

\DeclareMathOperator{\Id}{Id}

\DeclareMathOperator{\vect}{vect}

\renewcommand{\mod}{\operatorname{mod}}
\renewcommand{\Im}{\operatorname{Im}}

\renewcommand{\dim}{\operatorname{dim}}

\newcommand{\uvect}{\operatorname{\underline{vect}}}

\newcommand{\lto}[1][{}]{\stackrel{#1}{\longrightarrow}}
\newcommand{\rto}[1][{}]{\stackrel{#1}{\longleftarrow}}
\newcommand{\xto}{\xrightarrow}

\newcommand{\op}{\mathrm{op}}

\newcommand{\can}{\mathrm{can}}

\def\d{\delta}

\def\p{\phi}
\def\r{\rho}

\def\la{\lambda}

\def\Ga{\Gamma}

\def\Si{\Sigma}

\def\A{{\mathcal A}}
\def\B{{\mathcal B}}
\def\C{{\mathcal C}}

\def\Oc{{\mathcal O}}

\def\T{{\mathcal T}}

\def\bbP{\mathbb P}

\def\bbX{\mathbb X}

\def\bbZ{\mathbb Z}

\def\bfD{\mathbf D}

\def\bfL{\mathbf L}
\def\bfp{\mathbf p}

\newcommand{\bflambda}{\boldsymbol\lambda}

\newcommand{\ox}{\vec{x}}

\newcommand{\oc}{\vec{c}}

\begin{document}

\begin{abstract}
Expansions of abelian categories are introduced. These are certain
functors between abelian categories and provide a tool for
induction/reduction arguments. Expansions arise naturally in the study
of coherent sheaves on weighted projective lines; this is illustrated
by various applications.
\end{abstract}

\maketitle

\setcounter{tocdepth}{1}

\section{Introduction}\label{se:intro}

In this note we discuss certain functors between abelian categories
which we call \emph{expansions}. Such functors arise naturally in
finite dimensional representation theory and provide a tool for
induction/reduction arguments. We describe the formal properties of
these functors and give some applications.

Roughly speaking, an expansion is a fully faithful and exact functor
$\B\to\A$ between abelian categories that admits an exact left adjoint
and an exact right adjoint. In addition one requires the existence of
simple objects $S_\la$ and $S_\rho$ in $\A$ such that
$S_\la^\perp=\B={^\perp S_\rho}$, where $\B$ is viewed as a full
subcategory of $\A$ and
\begin{align*}
S_\la^\perp&=\{A\in\A\mid \Hom_\A(S_\la,A)=0=\Ext^1_\A(S_\la,A)\},\\
^\perp S_\r&=\{A\in\A\mid \Hom_\A(A,S_\r)=0=\Ext^1_\A(A,S_\r)\}.
\end{align*}
In fact, these simple objects are related by an
almost split sequence $0\to S_\rho\to S\to S_\la\to 0$ in $\A$ such that $S$
is a simple object in $\B$. In terms of the Ext-quivers of $\A$ and
$\B$, the expansion $\B\to\A$ turns the vertex $S$ into an arrow
$S_\la\to S_\rho$.

It is interesting to note that an expansion $\B\to\A$ induces a
recollement of derived categories
\[\xymatrixrowsep{3pc} \xymatrixcolsep{2pc}\xymatrix{
\bfD^b(\B)\;\ar[rr]|-{}&&\;\bfD^b(\A)\; \ar[rr]|-{}
\ar@<1.2ex>[ll]^-{}\ar@<-1.2ex>[ll]_-{}&&
\;\bfD^b(\mod\Delta)\ar@<1.2ex>[ll]^-{}\ar@<-1.2ex>[ll]_-{} }\] where
$\mod\Delta$ denotes the category of finitely generated right modules
over the \emph{associated division ring} $\Delta=\End_\A(S_\la)$ which
is isomorphic to $\End_\A(S_\rho)$.

Our motivation for studying expansions is the following result that
characterizes the abelian categories arising as categories of coherent
sheaves on weighted projective lines in the sense of Geigle and
Lenzing \cite{GL1987}.

\begin{Thm}
Let $k$ be an arbitrary field. A $k$-linear abelian
category $\A$ is equivalent to $\coh\bbX$ for some weighted projective
line $\bbX$ over $k$ (the exceptional points being rational with residue field $k$)
if and only if there exists a finite sequence
$\A^0\subseteq\A^1\subseteq \ldots\subseteq\A^r=\A$ of full
subcategories such that $\A^0$ is equivalent to $\coh\bbP^1_k$ and
each inclusion $\A^l\to\A^{l+1}$ is a non-split expansion with
associated division ring $k$.
\end{Thm}

This result is part of Theorem~\ref{th:char} and based on a technique
which is known as \emph{reduction of weights} \cite{GL1991}. There is
also an explicit construction for categories of coherent sheaves which
increases the weight type and is therefore called \emph{insertion of
  weights} \cite{Le1998}. The reduction of weights involves
\emph{perpendicular categories}; these were introduced in the context
of quiver representations by Schofield \cite{Sch} but appear already
in early work of Gabriel \cite{G}. A common feature of such reduction
techniques is the reduction of the rank of the Grothendieck
group. Note that in most cases the Grothendieck group is free of
finite rank.  Further induction/reduction techniques in representation
theory include \emph{one-point extensions of algebras} and the
\emph{shrinking of arrows}; see Ringel's report on tame algebras in
\cite{Ri1980}. In a more categorical setting, the \emph{trivial
  extensions} introduced by Fossum, Griffith and Reiten \cite{FGR}
should be mentioned. One-point extensions are probably the most
important among these techniques, but they require the existence of
enough projective objects. A category with Serre duality has no
non-zero projective objects, and in that sense expansions seem to be
the appropriate variant for dealing with sheaves on weighted
projective lines.

This paper is organized as follows. Some preliminaries about abelian
categories are collected in Section~\ref{se:prelim}. The central
subject of this work is treated in Section~\ref{se:expansions},
including a couple of generic examples. The applications to weighted
projective lines are discussed in the final Section~\ref{se:appl}.

\section{Preliminaries}\label{se:prelim}

In this section we fix our notation and collect some basic facts about
abelian categories. The standard reference is \cite{G}.

Let $F\colon \A\to\B$ be an additive functor between additive
categories. The \emph{kernel} $\Ker F$ of $F$ is the full
subcategory of $\A$ formed by all objects $A$ such that
$FA=0$. The \emph{essential image} $\Im F$ of $F$ is the full
subcategory of $\B$ formed by all objects $B$ such that $B$
is isomorphic to $FA$ for some $A$ in $\A$. Observe that $F$
is fully faithful if and only if $F$ induces an equivalence
$\A\xto{\sim}\Im F$; in this case we usually identify
$\A$ with $\Im F$, and identify $F$ with the inclusion
functor $\Im F\rightarrowtail \B$.

\subsection{Serre subcategories and quotient categories}
A non-empty full subcategory $\C$ of an abelian category $\A$ is a
\emph{Serre subcategory} provided that $\C$ is closed under taking
subobjects, quotients and extensions. This means that for every exact
sequence $0\to A'\to A\to A''\to 0$ in $\A$, the object $A$ belongs to
$\C$ if and only if $A'$ and $A''$ belong to $\C$.

Given a Serre subcategory $\C$ of $\A$, the \emph{quotient category}
$\A/\C$ is by definition the localization of $\A$ with respect to the
collection of morphisms that have their kernel and cokernel in
$\C$. The quotient category $\A/\C$ is an abelian category and the
quotient functor $Q\colon \A\to \A/\C$ is exact with $\Ker Q=\C$.

We observe that the kernel $\Ker F$ of an exact functor $F\colon \A\to
\B$ between abelian categories is a Serre subcategory of $\A$. Given
any Serre subcategory $\C\subseteq \Ker F$, the functor $F$ induces a
unique functor $\bar{F}\colon \A/\C\to \B$ such that $F=\bar{F}
Q$; moreover, the functor $\bar{F}$ is exact.

\subsection{Perpendicular categories}

Let $\A$ be an abelian category.  The quotient functor $\A\to \A/\C$
with respect to a Serre subcategory $\C$ admits an explicit
description if there exists a right adjoint; this is based on the use
of the perpendicular category $\C^\perp$.

For any class $\C$ of objects in $\A$, its \emph{perpendicular
  categories} are by definition the full subcategories
\begin{align*}
\C^\perp&=\{A\in\A\mid\Hom_\A(C,A)=0=\Ext^1_\A(C,A)\text{ for all
}C\in\C\},\\
^\perp\C&=\{A\in\A\mid\Hom_\A(A,C)=0=\Ext^1_\A(A,C)\text{ for all
}C\in\C\}.
\end{align*}

The next result shows that this definition of a perpendicular category
is appropriate in the abelian context. The lemma provides a useful
criterion for an exact functor to be a quotient functor and it
describes the right adjoint of a quotient functor.

\begin{lem}[{\cite[Chap.\ III.2]{G}}]\label{le:perp}
\pushQED{\qed} Let $F\colon\A\to\B$ be an exact functor between abelian categories
and suppose that $F$ admits a right adjoint $G\colon\B\to\A$. Then the
following are equivalent:
\begin{enumerate}
\item The functor $F$ induces an equivalence $\A/\Ker F\xto{\sim}\B$.
\item The functor $F$ induces an equivalence $(\Ker F)^\perp\xto{\sim}\B$.
\item The functor $G$ induces an equivalence $\B\xto{\sim}(\Ker F)^\perp$.
\item The functor $G$ is fully faithful.
\end{enumerate}
Moreover, in that case $(\Ker F)^\perp=\Im G$ and  $\Ker F={^\perp}(\Im G)$.\qedhere
\end{lem}

Next we characterize the fact that a quotient functor admits a right adjoint.

\begin{lem}[{\cite[Prop.\ 2.2]{GL1991}}]\label{le:section}
\pushQED{\qed} Let $\A$ be an abelian category and $\C$ a Serre
subcategory. Then the quotient functor $\A\to\A/\C$ admits a right
adjoint if and only if every object $A$ in $\A$ fits into
an exact sequence
\begin{equation}\label{eq:local}
0\lto A'\lto A\lto \bar A\lto A''\lto 0
\end{equation}
such that $A',A''\in\C$ and $\bar A\in\C^\perp$.

In that case the functor $\A\to\C$ sending $A$ to $A'$ is a
right adjoint of the inclusion $\C\to\A$, and
the functor $\A\to\C^\perp$ sending $A$ to $\bar A$ is a left adjoint of
the inclusion $\C^\perp\to\A$.\qedhere
\end{lem}

\subsection{Extensions}

Let $\A$ be an abelian category. For a pair of objects $A,B$ and
$n\geq 1$, let $\Ext^n_\A(A, B)$ denote the group of extensions in the
sense of Yoneda. Recall that an
element $[\xi]$ in $\Ext^n_\A(A, B)$ is represented by
an exact sequence $\xi\colon 0\rightarrow B\rightarrow
E_n\rightarrow \cdots \rightarrow E_1 \rightarrow A\rightarrow 0$ in
$\A$. Set $\Ext^0_\A(A, B)=\Hom_\A(A, B)$.

\begin{lem}\label{le:ext}
Let $F\colon\A\to\B$ and $G\colon\B\to\A$ be a pair of exact functors
such that $F$ is a left adjoint of $G$. Then we have natural
isomorphisms
\[\Ext^n_\B(FA, B)\cong \Ext_\A^n(A,GB)\]
for all $A\in \A$, $B\in\B$ and $n\geq 0$.
\end{lem}
\begin{proof}
The case $n=0$ is clear. For $n\geq 1$, the isomorphism sends $[\xi]$
in $\Ext^n_\B(FA,B)$ to $[(G\xi).\eta_A]$ in $\Ext_\A^n(A,GB)$, where
$\eta_A\colon A\rightarrow GF(A)$ is the unit of the adjoint pair and
$(G\xi).\eta_A$ denotes the pullback of $G\xi$ along $\eta_A$.
\end{proof}

\section{Expansions of abelian categories}\label{se:expansions}

In this section we introduce the concept of expansion and contraction
for abelian categories.\footnote{The authors are indebted to Claus
Michael Ringel for suggesting the terms `expansion' and
`contraction'.}  Roughly speaking, an expansion is a fully faithful
and exact functor $\B\to\A$ between abelian categories that admits an
exact left adjoint and an exact right adjoint. In addition one
requires the existence of simple objects $S_\la$ and $S_\rho$ in $\A$
such that $S_\la^\perp=\B={^\perp S_\rho}$, where $\B$ is viewed as a
full subcategory of $\A$. In fact, these simple objects are related by
an exact sequence $0\to S_\rho\to S\to S_\la\to 0$ in $\A$ such that
$S$ is a simple object in $\B$. In terms of the Ext-quivers of $\A$
and $\B$, the expansion $\B\to\A$ turns the vertex $S$ into an arrow
$S_\la\to S_\rho$. On the other hand, $\B$ is a contraction of $\A$ in
the sense that the left adjoint of $\B\to\A$ identifies $S_\rho$ with
$S$, whereas the right adjoint identifies $S_\la$ with $S$.

In the following we use the term `expansion' but there are interesting
situations where `contraction' yields a more appropriate point of
view. So one should think of a process having two directions that are
opposite to each other.

\subsection{Left and right expansions}

Let $\A$ be an abelian category.  A full subcategory $\B$ of $\A$ is
called \emph{exact abelian} if $\B$ is an abelian category and the
inclusion functor is exact. Thus a fully faithful and exact functor
$\B\to\A$ between abelian categories identifies $\B$ with an exact
abelian subcategory of $\A$.

\begin{defn}
A fully faithful and exact functor $\B\to\A$ between abelian
categories is called \emph{left expansion} if the following conditions
are satisfied:
\begin{enumerate}
\item[{\rm (E1)}] The functor $\B\to\A$ admits an exact left adjoint.
\item[{\rm (E2)}] The category $^\perp\B$ is equivalent to
  $\mod\Delta$ for some division ring $\Delta$.
\item[{\rm (E3)}] $\Ext^2_\A(A,B)=0$ for all $A,B\in{^\perp\B}$.
\end{enumerate}
The functor $\B\to\A$ is called  \emph{right
expansion} if the dual conditions are satisfied.
\end{defn}

\begin{lem}\label{le:leftmax}
Let $i\colon\B\to\A$ be a left expansion of abelian categories. Denote
by $i_\la$ its left adjoint and set $\C=\Ker i_\la$. Then the
following holds.
\begin{enumerate}
\item The category $\C$ is a Serre subcategory of $\A$ satisfying
$\C={^\perp\B}$ and $\C^\perp=\B$.
\item The composite $\B\xto{i}\A\xto{\can}\A/\C$ is
an equivalence and the left adjoint $i_\la$
induces a quasi-inverse $\A/\C\xto{\sim}\B$.
\item $\Ext^n_\B(i_\la A,B)\cong\Ext^n_\A(A,iB)$ for all $A\in\A$, $B\in\B$, and $n\geq 0$.
\end{enumerate}
\end{lem}
\begin{proof}
(1) and (2) follow from Lemma~\ref{le:perp}, while (3) follows
from Lemma~\ref{le:ext}.
\end{proof}

\begin{defn}
An object $S$ in an abelian category $\A$ is called \emph{localizable}
if the following conditions are satisfied:
\begin{enumerate}
\item[{\rm (L1)}] The object $S$ is simple.
\item[{\rm (L2)}] $\Hom_\A(S,A)$ and $\Ext^1_\A(S,A)$ are of finite length over
$\End_\A(S)$ for all $A\in\A$.
\item[{\rm (L3)}] $\Ext^1_\A(S,S)=0$ and $\Ext^2_\A(S,A)=0$ for all $A\in\A$.
\end{enumerate}
The object $S$ is called \emph{colocalizable} if the dual conditions
are satisfied.
\end{defn}

The following lemma describes for any abelian category $\A$ a
bijective correspondence between localizable objects in $\A$ and left
expansions $\B\to\A$.

For an object $X$ in $\A$, we denote by $\add X$ the full subcategory
consisting of all finite direct sums of copies of $X$ and their direct
summands.

\begin{lem}\label{le:localizable}
Let $\A$ be an abelian category.
\begin{enumerate}
\item If $i\colon\B\to\A$ is a left expansion, then there exists a
  localizable object $S\in\A$ such that $S^\perp=\Im i$.
\item If $S\in\A$ is a localizable object, then the inclusion
  $S^\perp\rightarrowtail\A$ is a left expansion.
\end{enumerate}
\end{lem}
\begin{proof}
(1) We identify $\B=\Im i$. Let $S$ be an indecomposable object in
  $^\perp\B$.  Then $S$ is a simple object and $\Ext_\A^1(S,S)=0$
  since $^\perp\B=\add S$ is semisimple.  For each object $A$ in $\A$,
  we use the natural exact sequence \eqref{eq:local}
$$0\lto A'\lto A\lto[\eta_A] \bar A\lto A''\lto 0$$ with $A',A''\in
{^\perp\B}$ and $\bar A\in\B$. This sequence induces the following
isomorphisms.
\begin{gather*}\Hom_\A(S,A')\lto[\sim]\Hom_\A(S,A)\\
\Ext_\A^1(S,A)\lto[\sim]\Ext_\A^1(S,\Im\eta_A)\rto[\sim]\Hom_\A(S,A'')
\end{gather*}
Here we use the condition $\Ext^2_\A(S,S)=0$.  It follows that
$\Hom_\A(S,A)$ and $\Ext^1_\A(S,A)$ are of finite length over
$\End_\A(S)$. Now observe that the functor sending $A$ to
$\Hom_\A(S,A'')$ is right exact. Thus $\Ext^1_\A(S,-)$ is right exact,
and therefore $\Ext^2_\A(S,-)=0$, for example by
\cite[Lemma~A.1]{RvB}.  Finally, $S^\perp=\B$ follows from
Lemma~\ref{le:perp}.

(2) The proof follows closely that of \cite[Prop.~3.2]{GL1991}. Set
$\B=S^\perp$ and observe that this is an exact abelian subcategory
since $\Ext_\A^2(S,-)=0$. A left adjoint $i_\la$ of the inclusion
$\B\rightarrowtail\A$ is constructed as follows.  Fix an object $A$ in
$\A$.  There exists an exact sequence $0\to A\to B\to S^n\to 0$ for
some $n\ge 0$ such that $\Ext_\A^1(S,B)=0$ since $\Ext^1_\A(S,A)$ is
of finite length over $\End_\A(S)$. Now choose a morphism $S^m\to B$
such that the induced map $\Hom_\A(S,S^m)\to\Hom_\A(S,B)$ is
surjective and let $\bar A$ be its cokernel. It is easily checked that
the composite $A\to B\to \bar A$ is the universal morphism into
$\B$. Thus we define $i_\la A=\bar A$.

Next observe that kernel and cokernel of the adjunction morphism $A\to
i_\la A$ belong to $\C=\add S$ for each object $A$ in $\A$. Moreover,
$\C$ is a Serre subcategory of $\A$ since $S$ is simple and
$\Ext^1_\A(S,S)=0$. Thus we can apply Lemma~\ref{le:section} and infer
that the quotient functor $\A\to\A/\C$ admits a right adjoint. In fact
the right adjoint identifies $\A/\C$ with $\C^\perp$, by
Lemma~\ref{le:perp}, and therefore $i_\la$ identifies with the
quotient functor. In particular, $i_\la$ is exact. We have
$^\perp\B=\C$ by Lemma~\ref{le:perp}, and $\Hom_\A(S,-)$ induces an
equivalence $\C\xto{\sim}\mod\End_\A(S)$.  Thus the inclusion
$\B\rightarrowtail\A$ is a left expansion.
\end{proof}

\subsection{Expansions of abelian categories}

We are ready to introduce the central concept of this work.

\begin{defn}
A fully faithful and exact functor between abelian
categories is called an \emph{expansion} of
abelian categories if the functor is a left and a right expansion.
\end{defn}

Let us fix some notation for an expansion $i\colon \B\to\A$.  We
identify $\B$ with the essential image of $i$. We denote by $i_\la$
the left adjoint of $i$ and by $i_\rho$ the right adjoint of $i$. We
choose an indecomposable object $S_\la$ in $^\perp\B$ and an
indecomposable object $S_\rho$ in $\B^\perp$. Thus $^\perp\B=\add
S_\la$ and $\B^\perp=\add S_\rho$. Finally, set $\bar
S=i_\la(S_\rho)$.

\begin{prop}\label{pr:localizable}
Let $\A$ be an abelian category.
\begin{enumerate}
\item Given an expansion $i\colon\B\to\A$, there exist a localizable
  object $S_\la$ and a colocalizable object $S_\r$ such that
  $S_\la^\perp=\Im i={^\perp S_\r}$.
\item Let $S_\la$ be a localizable object and $S_\r$ a colocalizable
  object in $\A$ such that $S_\la^\perp={^\perp S_\r}$.  Then the
  inclusion  $S_\la^\perp\rightarrowtail\A$ is an expansion.
\end{enumerate}
\end{prop}

\begin{proof} Apply
Lemma~\ref{le:localizable} and its dual.
\end{proof}

An expansion $i\colon\B\to\A$ is called \emph{split} if
$\B^\perp={^\perp\B}$. If the expansion is non-split, then the exact
sequences \eqref{eq:local} for $S_\la$ and $S_\rho$ are of the form
\begin{equation}\label{eq:locseq}
0\to S_\rho\to ii_\la(S_\rho)\to S_\la^l\to 0\quad\text{and}\quad
0\to S_\rho^r\to ii_\rho(S_\la)\to S_\la\to 0
\end{equation}
for some integers $l,r\geq 1$. In Lemma~\ref{le:extsimple}, we see that $l=1=r$.

\begin{lem}
Let $\B\to\A$ be an expansion of abelian categories. Then the
following are equivalent:
\begin{enumerate}
\item The expansion $\B\to\A$ is split.
\item $\A=\B\amalg\C$ for some Serre subcategory $\C$ of $\A$.
\item $\B$ is a Serre subcategory of $\A$.
\end{enumerate}
\end{lem}
\begin{proof}
(1) $\Rightarrow$ (2): Take $\C={^\perp\B}=\B^\perp$.

(2) $\Rightarrow$ (3): An object $A\in\A$ belongs to $\B$ if and only
if $\Hom_\A(A,B)=0$ for all $B\in\C$. Thus $\B$ is closed under taking
quotients and extensions. The dual argument shows that $\B$ is closed
under taking subobjects.

(3) $\Rightarrow$ (1): If the expansion is non-split, then the
sequences in \eqref{eq:locseq} show that $\B$ is not a Serre
subcategory.
\end{proof}

\begin{lem}\label{le:extsimple}
Let $\B\to\A$ be a non-split expansion of abelian categories.
\begin{enumerate}
\item The object $\bar S=i_\la(S_\rho)$ is a simple object in $\B$ and
isomorphic to $i_\rho(S_\la)$.
\item The functor $i_\la$ induces an equivalence
 $\B^\perp\xto{\sim}\add\bar S$.
\item The functor $i_\rho$ induces an equivalence
$^\perp\B\xto{\sim}\add\bar S$.
\end{enumerate}
\end{lem}
\begin{proof}
(1) Let $\p\colon i_\la(S_\rho)\to A$ be a non-zero morphism in
$\B$. Adjunction takes this to a monomorphism $S_\rho\to A$ in
$\A$ since $S_\rho$ is simple. Applying $i_\la$ gives back a morphism which is isomorphic to
$\p$. This is a monomorphism since $i_\la$ is exact. Thus
$i_\la(S_\rho)$ is simple.

Now apply $i_\rho$ to the first and $i_\la$ to the second sequence in
\eqref{eq:locseq}. Then
$$i_\la(S_\rho)\cong i_\rho(S_\la)^l\cong i_\la(S_\rho)^{lr}.$$
This implies $l=1=r$ and therefore $i_\la(S_\rho)\cong i_\rho(S_\la)$.

(2) We have a sequence of isomorphisms
$$\Hom_\A(S_\rho,S_\rho)\lto[\sim]\Hom_\A(S_\rho,ii_\la(S_\rho))\lto[\sim]
\Hom_\B(i_\la(S_\rho),i_\la(S_\rho))$$ which takes a morphism $\p$ to
$i_\la\p$. Thus $i_\la$ induces an equivalence $\add
S_\rho\xto{\sim}\add i_\la(S_\rho)$.

(3) Follows from (2) by duality.
\end{proof}

An expansion $i\colon\B\to\A$ of abelian categories determines a
division ring $\Delta$ such that $^\perp\B$ and $\B^\perp$ are
equivalent to $\mod\Delta$; we call $\Delta$ the \emph{associated
  division ring}.  If the expansion does not split, then the sequences
in \eqref{eq:locseq} yield an essentially unique non-split
extension \[0\lto S_\r\to i\bar S\lto S_\la\lto 0\] which is called
the \emph{connecting sequence} of the expansion $i$.  This sequence is
almost split; see Proposition~\ref{pr:ARformula}.

\subsection{Recollements}

Fix an expansion $i\colon\B\to\A$ with associated division ring
$\Delta$. We identify the perpendicular categories of $\B$ with
$\mod\Delta$ via the equivalences
${^\perp\B}\xto{\sim}\mod\Delta\xleftarrow{\sim}\B^\perp$. There are
inclusions $j\colon{^\perp\B}\to\A$ and $k \colon{\B^\perp}\to\A$ with
adjoints $j_\rho$ and $k_\la$. These functors yield the following
diagram.
\[\xymatrixrowsep{3pc} \xymatrixcolsep{2pc}\xymatrix{
\B\;\ar[rr]|-{i}&&\;\A\;
\ar@<1.0ex>[rr]|-{j_\rho}\ar@<-1.0ex>[rr]|-{k_\la}
\ar@<1.2ex>[ll]^-{i_\rho}\ar@<-1.2ex>[ll]_-{i_\la}&&
\;\mod\Delta\ar@<2.2ex>[ll]^-{k}\ar@<-2.2ex>[ll]_-{j} }\] It is
interesting to observe that this diagram induces a recollement of
triangulated categories \cite{BBD} when one passes from abelian categories
to their derived categories.

Recall that a diagram of exact functors between triangulated
categories
\begin{equation*}\label{eq:recollement}
\xymatrixrowsep{3pc} \xymatrixcolsep{2pc}\xymatrix{
  \mathcal{T}'\;\ar[rr]|-{i}&&\;\mathcal{T}\; \ar[rr]|-{j}
  \ar@<1.2ex>[ll]^-{i_\r}\ar@<-1.2ex>[ll]_-{i_\la}&&
  \;\mathcal{T}''\ar@<1.2ex>[ll]^-{j_\r}\ar@<-1.2ex>[ll]_-{j_\la}
}\end{equation*} forms a \emph{recollement}, provided that the
following conditions are satisfied:
\begin{enumerate}
\item[{\rm (R1)}] The pairs $(i_\la, i)$, $(i, i_\rho)$, $(j_\la, j)$, and $(j,j_\rho)$ are adjoint.
\item[{\rm (R2)}] The functors $i$, $j_\la$, and $j_\rho$ are fully faithful.
\item[{\rm (R3)}] $\Im i=\Ker j$.
\end{enumerate}
Note that in this case $j$ is a quotient functor in the sense of
Verdier, inducing an equivalence $\T/\Ker j\xto{\sim}\T''$.

Given any abelian category $\A$, we denote by $\bfD^b(\A)$ its bounded
derived category. An exact functor $F\colon\A\to\B$ between abelian
categories extends to an exact functor $\bfD^b(\A)\to\bfD^b(\B)$ which
we denote by $F^*$.  Note that kernel of $F^*$ coincides with the full
subcategory $\bfD^b_{\Ker F}(\A)$ consisting of the complexes in
$\bfD^b(\A)$ with cohomology in $\Ker F$.

The following lemma describes the functors that yield a recollement of
derived categories. Compare the first part with Lemma~\ref{le:ext}.

\begin{lem}\label{le:der-adjoint}
Let $F\colon \B\to\A$ be a fully faithful exact functor between
abelian categories and suppose that $F$ admits an exact right adjoint
$G\colon\A\to\B$.
\begin{enumerate}
\item $F^*$ and $G^*$ form an adjoint pair of exact functors and
  $F^*$ is fully faithful.
\item The inclusion $\Ker G^*\to\bfD^b(\A)$ admits a left adjoint
  which induces an equivalence $\bfD^b(\A)/\Im F^*\xto{\sim}\Ker G^*$.
\end{enumerate}
\end{lem}
\begin{proof}
(1) Fix a pair of complexes $X\in\bfD^b(\B)$ and $Y\in\bfD^b(\A)$. The
  unit $\Id_\B\to GF$ yields a morphism $\eta_X\colon X\to G^*F^*(X)$
  and we obtain a map
\[\Hom_{\bfD^b(\A)}(F^*X,Y)\lto\Hom_{\bfD^b(\B)}(X,G^*Y)\]
by sending $\p$ to $(G^*\p)\eta_X$. A simple application of Beilinson's
lemma shows that this map is bijective. The same lemma yields that
$F^*$ is fully faithful.

(2) We construct a left adjoint $L\colon \bfD^b(\A)\to \Ker G^*$ as
follows. For $X\in\bfD^b(\A)$ complete the counit $F^*G^*(X)\to X$ to
an exact triangle $F^*G^*(X)\to X\to X'\to$. The assignment $X\mapsto
L(X)=X'$ is functorial and yields an exact left adjoint for the
inclusion $\Ker G^*\to\bfD^b(\A)$; see for example
\cite[Lemma~3.3]{BIK} for details. A right adjoint of a fully faithful
functor is, up to an equivalence, a quotient functor, by
\cite[Prop.~I.1.3]{GZ}. Thus it remains to observe that $\Ker L=\Im
F^*$, which is obvious from the triangle defining $L$.
\end{proof}

\begin{prop}\label{pr:recoll}
An expansion of abelian categories $i\colon\B\to\A$ with associated
division ring $\Delta$ induces the following recollement.
\[\xymatrixrowsep{3pc} \xymatrixcolsep{2pc}\xymatrix{
\bfD^b(\B)\;\ar[rr]|-{i^*}&&\;\bfD^b(\A)\; \ar[rr]|-{}
\ar@<1.2ex>[ll]^-{(i_\rho)^*}\ar@<-1.2ex>[ll]_-{(i_\la)^*}&&
\;\bfD^b(\mod\Delta)\ar@<1.2ex>[ll]^-{k^*}\ar@<-1.2ex>[ll]_-{j^*}
}\]
\end{prop}

We point out that in general the unlabeled functor is not induced from
an exact functor between the abelian categories. In fact, this functor
equals the right derived functor of $j_\r$, and it equals the left
derived functor of $k_\la$.

\begin{proof}
The assertion is an immediate consequence of
Lemma~\ref{le:der-adjoint} and its dual. It only remains to observe
that the inclusion $j$ induces an equivalence
$\bfD^b(\mod\Delta)\xto{\sim}\bfD^b_{\Ker i_\la}(\A)$, while $k$
induces an equivalence $\bfD^b(\mod\Delta)\xto{\sim}\bfD^b_{\Ker
  i_\r}(\A)$. This follows from an application of Beilinson's lemma.
\end{proof}

\subsection{Simple objects}
Let $i\colon\B\to\A$ be an expansion. The left adjoint $i_\la$ induces
a bijection between the isomorphism classes of simple objects of $\A$
that are different from $S_\la$, and the isomorphism classes of simple
objects of $\B$.  On the other hand, all simple objects of $\A$
correspond to simple objects of $\B$ via $i$.  All this is made
precise in the next lemma.

\begin{lem}\label{le:simple}
Let $i\colon\B\to\A$ be an expansion of abelian categories.
\begin{enumerate}
\item If $S$ is a simple object in $\B$ and $S\not\cong \bar S$, then
$iS$ is simple in $\A$ and $i_\la iS\cong S$.
\item There is an exact sequence $0\to S_\rho\to i\bar S\to S_\la\to
0$ in $\A$, provided the expansion $\B\to\A$ is non-split.
\item If $S$ is a simple object in $\A$ and $S\not\cong S_\la$, then
$i_\la S$ is simple in $\B$. Moreover, $S\cong i i_\la S$ if
$S\not\cong S_\rho$.
\end{enumerate}
\end{lem}
\begin{proof}
(1) Let $0\neq U\subseteq i S$ be a subobject. Then $\Hom_\B(i_\la
    U,S)\cong\Hom_\A(U,iS)\neq 0$ shows that $U\not\in\Ker i_\la$.
    Thus $i_\la U=S$, and therefore $iS/U$ belongs to $\Ker i_\la=\add
    S_\la$. On the other hand, $\Hom_\A(iS,S_\la)\cong\Hom_\B(S,\bar
    S)=0$. Thus $iS/U=0$, and it follows that $iS$ is simple. Finally
    observe that $i_\la iA\cong A$ for every object $A$ in $\B$.

(2) Take the exact sequence in \eqref{eq:locseq}.

(3) This is a general fact: A quotient functor $\A\to\A/\C$ takes each
    simple object of $\A$ not belonging to $\C$ to a simple object of
    $\A/\C$.  Here, we take $\C=\Ker i_\la$ and identify $i_\la$ with
    the corresponding quotient functor.

If $S\not\cong S_\rho$, then $i_\la S\not\cong \bar S$ and therefore
$ii_\la S$ is simple by (1). Thus the canonical morphism $S\to ii_\la
S$ is an isomorphism.
\end{proof}

The Ext-groups of most simple objects in $\A$ can be computed from
appropriate Ext-groups in $\B$. This follows from an adjunction
formula; see Lemma~\ref{le:leftmax}. The remaining cases are treated
in the following lemma.

\begin{lem}\label{le:ext_simple}
Let $i\colon\B\to\A$ be a non-split expansion of abelian categories.
\begin{enumerate}
\item $\Hom_\A(S_\la,S_\la)\cong\Ext_\A^1(S_\la,S_\rho)\cong\Hom_\A(S_\rho,S_\rho)$.
\item $\Ext_\B^n(\bar S,\bar S)\cong\Ext_\A^n(S_\rho,S_\la)$ for $n\geq 1$.
\end{enumerate}
\end{lem}
\begin{proof}
(1) Applying $\Hom_\A(S_\la,-)$ to the first sequence in
\eqref{eq:locseq} yields the isomorphism
$\Hom_\A(S_\la,S_\la)\cong\Ext_\A^1(S_\la,S_\rho)$. The other
isomorphism is dual.

(2) We have
\[\Ext_\B^n(i_\la(S_\rho),i_\la(S_\rho))\cong\Ext_\A^n(S_\rho,ii_\la(S_\rho))
\cong\Ext_\A^n(S_\rho,S_\la),\] where the first isomorphism follows
from Lemma~\ref{le:leftmax} and the second from the first sequence in
\eqref{eq:locseq}.
\end{proof}

For an abelian category $\A$ denote by $\A_0$ the full subcategory
formed by all finite length objects; it is a Serre subcategory.

\begin{prop}\label{pr:length}
Let $i\colon\B\to\A$ be an expansion of abelian
categories.
\begin{enumerate}
\item The functor $i$ and its adjoints $i_\la$ and $i_\rho$
send  finite length objects to finite length objects.
\item The restriction $\B_0\to\A_0$ is an expansion of
abelian categories.
\item The induced functor $\B/\B_0\to\A/\A_0$ is an equivalence.
\end{enumerate}
\end{prop}
\begin{proof}
(1) follows from Lemma~\ref{le:simple} and (2) is an immediate
    consequence.

(3) Let $\C=\Ker i_\la$. The functor $i_\la$ induces an equivalence
    $\A/\C\xto{\sim}\B$. Moreover, $\C\subseteq\A_0$ and $i_\la$
    identifies $\A_0/\C$ with $\B_0$ by (1). Then it follows from a
    Noether isomorphism that $i_\la$ induces an equivalence
    $\A/\A_0\xto{\sim}\B/\B_0$. This is a quasi-inverse for the
    functor $\B/\B_0\to\A/\A_0$ induced by $i$.
\end{proof}

Let $\A$ be an abelian category. We call an object $A$
\emph{torsion-free} if $\Hom_\A(S, A)=0$ for each simple object
$S$. An object $A$ is called \emph{$1$-simple} if it becomes a simple
object in the quotient category $\A/{\A_0}$, equivalently, if for each subobject
$A'\subseteq A$ either $A'$ or $A/A'$ is of finite length, but not
both.

\begin{prop}\label{pr:vec}
Let $i\colon \B\to\A$ be a non-split expansion. Then the functor $i$ and its
adjoints $i_\la$ and $i_\r$ preserve torsion-free objects and $1$-simple objects.
\end{prop}

\begin{proof}
Using adjunctions, it is clear that $i$ and $i_\r$ preserve
torsion-free objects. Now fix a torsion-free object $A\in\A$ and a morphism
$\p\colon S\to i_\la A$ with $S\in\B$ simple. Consider the natural
exact sequence \eqref{eq:local} \[0\lto A'\lto A\lto ii_\la(A)\lto
A''\lto 0\] and observe that $A'=0$. The composite $iS\to ii_\la(A)\to
A''$ has a non-zero kernel, since $A''$ belongs to $\add S_\la$, and
 $i\p$ maps this kernel to $A$. This implies $\p=0$.

The statement on $1$-simple objects follows from
Proposition~\ref{pr:length}. Here we note that $i_\la$ and $i_\r$
induce functors $\A/\A_0\rightarrow \B/\B_0$ that are quasi-inverse to
the functor $\B/\B_0\rightarrow \A/\A_0$ induced by $i$.
\end{proof}

\subsection{The Ext-quiver}

The \emph{Ext-quiver} of an abelian category $\A$ is a valued quiver
$\Si=\Si(\A)$ which is defined as follows.  The set $\Si_0$ of
vertices is a fixed set of representatives of the isomorphism classes
of simple objects in $\A$. For a simple object $S$, let $\Delta(S)$
denote its endomorphism ring, which is a division ring.  There is an
arrow $S\to T$ with valuation $\d_{S,T}=(s,t)$ in $\Si$ if
$\Ext_\A^1(S,T)\neq 0$ with $s=\dim_{\Delta(S)}\Ext_\A^1(S,T)$ and
$t=\dim_{\Delta(T)^\op}\Ext_\A^1(S,T)$. We write $\d_{S,T}=(0,0)$ if
$\Ext_\A^1(S,T)=0$.

Given an expansion $\B\to\A$, the Ext-quiver $\Si(\A)$ can be
computed explicitly from the Ext-quiver $\Si(\B)$, and vice versa. The
following statement makes this precise.

\begin{prop}
Let $i\colon\B\to\A$ be a non-split expansion of abelian
categories. The functor induces a bijection
$$\Si_0(\B)\smallsetminus\{\bar
S\}\lto[\sim]\Si_0(\A)\smallsetminus\{S_\la,S_\rho\},$$ and for each
pair  $U,V\in\Si_0(\B)\smallsetminus\{\bar S\}$ the
following identities:
\[
\d_{iU,iV}=\d_{U,V},\quad \d_{iU,S_\la}=\d_{U,\bar S},\quad
\d_{S_\rho,iV}=\d_{\bar S,V}, \quad \d_{S_\rho,S_\la}=\d_{\bar S,\bar
  S},\quad \d_{S_\la,S_\rho}=(1,1).
\]
\end{prop}
\begin{proof}
Combine Lemmas~\ref{le:leftmax}, \ref{le:simple}, and
\ref{le:ext_simple}.
\end{proof}
The following diagram shows how $\Si(\B)$ and $\Si(\A)$ are related.
\setlength{\unitlength}{1.5pt}
\[
\begin{picture}(120,25)
\put(-55,8){$\Si(\B)$}
\put(-05,8){$\bar S$}
\put(-30,02){\circle*{1}}
\put(-30,10){\circle*{1}}
\put(-30,18){\circle*{1}}
\put(25,02){\circle*{1}}
\put(25,10){\circle*{1}}
\put(25,18){\circle*{1}}
\put(-22,21){\vector(2,-1){15}}
\put(-22,-01){\vector(2,1){15}}
\put(-22,10){\vector(1,0){15}}
\put(02,11){\vector(2,1){15}}
\put(02,09){\vector(2,-1){15}}
\put(55,8){$\Si(\A)$}
\put(105,8){$S_\la$}
\put(137,8){$S_\rho$}
\put(80,02){\circle*{1}}
\put(80,10){\circle*{1}}
\put(80,18){\circle*{1}}
\put(170,02){\circle*{1}}
\put(170,10){\circle*{1}}
\put(170,18){\circle*{1}}
\put(88,21){\vector(2,-1){15}}
\put(88,-01){\vector(2,1){15}}
\put(88,10){\vector(1,0){15}}
\put(115,10){\vector(1,0){20}}
\put(147,11){\vector(2,1){15}}
\put(147,09){\vector(2,-1){15}}
\put(118,13){$\scriptstyle{(1,1)}$}
\end{picture}
\]

\subsection{Examples}
We discuss two examples. The first one arises from the study of
coherent sheaves on weighted projective lines, while the second one
describes expansions for representations of quivers.

\begin{exm}\label{ex:wpl}
Let $k$ be a field and $\A$  a $k$-linear abelian
category with finite dimensional Hom and Ext spaces. Assume that
$\A$ has \emph{Serre duality}, that is, there is an
equivalence $\tau\colon \A\xto{\sim} \A$ with
functorial $k$-linear isomorphisms
$$D\Ext^1_\A(A, B)\cong \Hom_\A(B, \tau A)$$ for all $A, B$ in
$\A$, where $D=\Hom_k(-, k)$ denotes the standard $k$-duality. Let
$S_\lambda$ be a simple object in $\A$ with $\Ext_\A^1(S_\la,
S_\la)=0$ and set $S_\rho=\tau S_\la$. Then by Serre duality
$S_\la^\perp={^\perp S_\r}$. Note that the category $\A$ is
\emph{hereditary}, that is, $\Ext_\A^2(-,-)=0$. It follows that $S_\la$
is localizable and $S_\r$ is colocalizable. By
Proposition~\ref{pr:localizable} the inclusion
$S_\la^\perp\rightarrowtail \A$ is an expansion, and this is non-split
since $S_\la \ncong S_\r$.

A specific example of an abelian category having the above properties
is the category of finite dimensional nilpotent representations of a
quiver $\Ga_n$ of type $\tilde A_n$ with cyclic orientation. In that case,
$S_\la^\perp$ is equivalent to the category of finite dimensional
nilpotent representations of $\Ga_{n-1}$.
\end{exm}

\begin{exm}
Let $k$ be a field. Consider a finite quiver $\Ga$ having two vertices
$a,b$ that are joined by an arrow $\xi\colon a\to b$ which is
the unique arrow starting at $a$ and the unique arrow
terminating at $b$.  \setlength{\unitlength}{1.5pt}
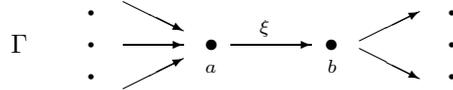
\[
\begin{picture}(80,25)
\put(-20,8){$\Ga$} \put(28.5,03){$\scriptstyle{a}$}
\put(59,03){$\scriptstyle{b}$} \put(30,10){\circle*{2.5}}
\put(60,10){\circle*{2.5}} \put(00,02){\circle*{1}}
\put(00,10){\circle*{1}} \put(00,18){\circle*{1}}
\put(90,02){\circle*{1}} \put(90,10){\circle*{1}}
\put(90,18){\circle*{1}} \put(08,21){\vector(2,-1){15}}
\put(08,-01){\vector(2,1){15}} \put(08,10){\vector(1,0){15}}
\put(35,10){\vector(1,0){20}} \put(67,11){\vector(2,1){15}}
\put(67,09){\vector(2,-1){15}} \put(42,13){$\scriptstyle{\xi}$}
\end{picture}
\]
We obtain a new quiver $\Ga'$ by identifying $a$ and $b$
and removing $\xi$.

Let $k\Ga$ be the path algebra of $\Ga$ and let $A=k\Ga/I$ be a finite
dimensional quotient algebra with respect to some admissible ideal
$I$. Denote by $\A=\mod A^\op$ the abelian category of finite
dimensional left $A$-modules, which is viewed as a full subcategory of
the category of $k$-linear representations of $\Ga$.

Consider the full subcategory $\B$ of $\A$ formed by all modules which
correspond to representations of $\Ga$ such that $\xi$ is represented
by an isomorphism. Note that $\B$ is equivalent to the category of
finite dimensional left modules over a finite dimensional algebra
$A'=k\Ga'/I'$ for some ideal $I'$.

For each vertex $i$ of $\Ga$, denote by $P_i$ (respectively, $I_i$)
the projective cover (respectively, injective hull) in $\A$ of the
corresponding simple module $S_i$. Note that the arrow $\xi$ induces
morphisms $P_\xi\colon P_b\to P_a$ and $I_\xi\colon I_b \to I_a$.

Assume that the morphism $P_\xi$ is a monomorphism and that $I_\xi$ is
an epimorphism. For example, this happens when the admissible ideal
$I$ is generated by paths which neither begin with nor end with
$\xi$.  Then we claim that the inclusion $\B\rightarrowtail\A$ is a
non-split expansion.

The assumption yields the following two exact sequences
\[0\to P_b\to P_a \to S_a\to
0 \quad\text{and}\quad 0\to S_b \to I_b\to I_a \to 0.\] It follows
that the simple module $S_a$ is a localizable object and that $S_b$ is
a colocalizible object of $\A$. Moreover, we observe from the two
exact sequences above that $S_a^\perp=\B={^\perp S_b}$. Thus by
Proposition~\ref{pr:localizable} the inclusion functor
$\B\rightarrowtail \A$ is a non-split expansion.

Observe that the algebra $A'$ is Morita equivalent to the universal
localization \cite{Sch1985} of $A$ at the map $P_\xi$, since $\B$
equals the full subcategory consisting of all $A$-modules $X$
such that $\Hom_A(P_\xi,X)$ is invertible.

Now let $A=k\Ga/I$ be an arbitrary finite dimensional algebra over
$k$, where $\Ga$ is any finite quiver and $I$ is an admissible ideal.
Then one can show that for any non-split expansion $i\colon\B\to
\A=\mod A^\op$, there is a unique arrow $\xi$ of $\Ga$ satisfying the
above conditions such that $i$ identifies $\B$ with the full
subcategory formed by all representations inverting $\xi$.
\end{exm}

\subsection{An Auslander-Reiten formula}
Given a non-split expansion $\B\to\A$, the corresponding
simple objects $S_\la$ and $S_\rho$ in $\A$ are related by an
Auslander-Reiten formula.

\begin{prop}\label{pr:ARformula}
Let $\B\to\A$ be a non-split expansion of abelian
categories and $\Delta$ its associated division ring. Then
$$D\Ext_\A^1(-,S_\rho)\cong\Hom_\A(S_\la,-)\quad\text{and}\quad
D\Ext_\A^1(S_\la,-)\cong\Hom_\A(-,S_\rho),$$ where
$D=\Hom_\Delta(-,\Delta)$ denotes the standard duality. In particular, any
non-split extension $0\to S_\rho\to E\to S_\la\to$ is an almost split
sequence.
\end{prop}
\begin{proof}
Recall that ${^\perp\B}=\add S_\la$ and $\B^\perp=\add S_\rho$.  Fix
an object $A$ in $\A$ and consider the corresponding exact sequence
\eqref{eq:local}
$$0\lto A'\lto A\lto \bar A\lto A''\lto 0.$$ with $A',A''$ in
${^\perp\B}$ and $\bar A$ in $\B$.  The morphism $A'\to A$ induces the
first and the third isomorphism in the sequence below, while the second
isomorphism follows from the isomorphism
$\Hom_\A(S_\la,S_\la)\cong\Ext_\A^1(S_\la,S_\rho)$ in
Lemma~\ref{le:ext_simple}.
$$D\Ext_\A^1(A,S_\rho)\cong
D\Ext_\A^1(A',S_\rho)\cong\Hom_\A(S_\la,A')\cong\Hom_\A(S_\la,A).$$
The isomorphism $D\Ext_\A^1(S_\la,-)\cong\Hom_\A(-,S_\rho)$ follows
from the first by duality.

Let $\xi\colon 0\to S_\rho\to E\to S_\la\to$ be a non-split extension.
Using the formula for $D\Ext^1_\A(-,S_\r)$, one shows that the
pullback along any morphism $\p\colon A\to S_\la$ is a split
extension, provided that $\p$ is not a split epimorphism. Analogously,
one shows via the formula for $D\Ext^1_\A(S_\la,-)$ that the pushout
along any morphism $\p\colon S_\r\to A$ is a split extension, provided
that $\p$ is not a split monomorphism. Thus $\xi$ is an almost split
sequence.
\end{proof}

\section{Applications to coherent sheaves on weighted projective lines}\label{se:appl}

In this section we present some results on weighted projective lines
that are based on expansions.  We begin with a brief description of
weighted projective lines and their categories of coherent sheaves.

Throughout this section we fix an arbitrary field $k$.

\subsection{Coherent sheaves on weighted projective lines}\label{se:coh}

Let $\bbP^1_k$ be the projective line over $k$, let
$\bflambda=(\la_1,\dots,\la_n)$ be a (possibly empty) collection of
distinct rational points of $\bbP^1_k$, and let $\bfp=(p_1,\dots,p_n)$
be a \index{weight sequence} \emph{weight sequence}, that is, a
sequence of positive integers.  The triple $\bbX =
(\bbP^1_k,\bflambda,\bfp)$ is called a \index{weighted projective
  line} \emph{weighted projective line}.  

Geigle and Lenzing \cite{GL1987} have associated to each weighted
projective line a category $\coh\bbX$ of coherent sheaves on $\bbX$,
which is the quotient category of the category of finitely generated
$\bfL(\bfp)$-graded $S(\bfp,\bflambda)$-modules, modulo the Serre
subcategory of finite length modules. Thus
\[\coh\bbX=\mod^{\bfL(\bfp)} S(\bfp,\bflambda)/(\mod^{\bfL(\bfp)}S(\bfp,\bflambda))_0.\]
Here $\bfL(\bfp)$ is the rank 1
additive group
\[
\bfL(\bfp) = \langle \ox_1,\dots,\ox_n,\oc \mid p_1 \ox_1 =
\dots = p_n \ox_n = \oc \rangle,
\]
and
\[
S(\bfp,\bflambda) = k[u,v,x_1,\dots,x_n] / ( x_i^{p_i} + \la_{i1} u -
\la_{i0} v ),
\]
with grading $\deg u = \deg v = \oc$ and $\deg x_i = \ox_i$, where
$\la_i = [\la_{i0}:\la_{i1}]$ in $\bbP^1_k$.  Geigle and Lenzing
showed that $\coh\bbX$ is a hereditary abelian category with finite
dimensional Hom and Ext spaces.  The free module $S(\bfp,\bflambda)$
yields a structure sheaf $\Oc$, and shifting the grading gives twists
$E(\ox)$ for any sheaf $E$ and $\ox\in\bfL(\bfp)$.

Every sheaf is the direct sum of a torsion-free sheaf and a finite
length sheaf. A torsion-free sheaf has a finite filtration by
\emph{line bundles}, that is, sheaves of the form $\Oc(\ox)$. The
finite length sheaves are easily described as follows.  There are
simple sheaves $S_x$ in bijection to closed points $x$ in
$\bbP^1_k\smallsetminus \bflambda$, and $S_{ij}$ ($1\le i\le n$, $1\le
j\le p_i$) satisfying for any $r\in\bbZ$ that
$\Hom(\Oc(r\oc),S_{ij})\neq 0$ if and only if $j=1$, and the only
extensions between them are
\[
\Ext^1(S_x,S_x) =k(x),
\quad
\Ext^1(S_{ij},S_{ij'}) = k
\quad
(j' \equiv j-1 \, (\mod \, p_i)).
\]
Here $k(x)$ denotes the residue field at each closed point $x$. 
For each simple sheaf $S$ and $l>0$ there is a unique sheaf with
length $l$ and top $S$, which is \emph{uniserial}, meaning that it has
a unique composition series.  These are all the finite length
indecomposable sheaves.

\subsection{A characterization of coherent sheaves on weighted projective lines}

The following result describes in terms of expansions the abelian
categories that arise as categories of coherent sheaves on weighted
projective lines.

\begin{thm}\label{th:char}
A $k$-linear abelian category $\A$ is equivalent to $\coh\bbX$ for
some weighted projective line $\bbX$ over $k$ if
and only if there exists a finite sequence $\A^0\subseteq\A^1\subseteq
\ldots\subseteq\A^r=\A$ of full subcategories such that $\A^0$ is
equivalent to $\coh\bbP^1_k$ and each inclusion $\A^l\to\A^{l+1}$ is a
non-split expansion with associated division ring $k$.

In that case each inclusion $\A^l\to\A^{l+1}$ induces a recollement
\[\xymatrixrowsep{3pc} \xymatrixcolsep{2pc}\xymatrix{
\bfD^b(\A^l)\;\ar[rr]|-{}&&\;\bfD^b(\A^{l+1})\; \ar[rr]|-{}
\ar@<1.2ex>[ll]^-{}\ar@<-1.2ex>[ll]_-{}&&
\;\bfD^b(\mod k).\ar@<1.2ex>[ll]^-{}\ar@<-1.2ex>[ll]_-{}
}\]
\end{thm}

\begin{proof}
The first part of the assertion is covered by
\cite[Thm.~6.8.1]{CK} and based on work of Lenzing in \cite{L}. A
detailed proof can be found in \cite{CK}. The existence of a
recollement induced by the inclusion $\A^l\to\A^{l+1}$ is an immediate
consequence of Proposition~\ref{pr:recoll}.

Now fix a weighted projective line $\bbX = (\bbP^1_k,\bflambda,\bfp)$.
We provide the argument that gives the filtration
$\A^0\subseteq\A^1\subseteq \ldots\subseteq\A^r=\A$ for $\A=\coh\bbX$;
this is known as \emph{reduction of weights}. If $p_i=1$ for all $i$,
then $\A=\coh\bbP^1_k$. Otherwise, choose some $j$ such that
$p_j>1$. Then $S=S_{j1}$ is a simple object satisfying
$\Ext^1_\A(S,S)=0$. So we can apply Example~\ref{ex:wpl} and obtain an
expansion $S^\perp\rightarrowtail\A$. It follows from the arguments
given in \cite[Prop.~1]{L} that $S^\perp$ is equivalent to $\coh\bbX'$
for $\bbX' = (\bbP^1_k,\bflambda,\bfp')$, where $p_i'=p_i-\delta_{ij}$
for $1\le i\le n$. Thus we can proceed and obtain a trivial weight
sequence $(1,\ldots,1)$ after $\sum_{i=1}^n(p_i-1)$ steps.

Next we provide an explicit construction of an expansion
$\coh\bbX'\to\coh\bbX$, where $\bbX' = (\bbP^1_k,\bflambda,\bfp')$
with $p_i'=p_i-\delta_{ij}$ for $1\le i\le n$ and some fixed $j$; see
also \cite[\S9]{GL1991}.  There is an injective map
$$\phi\colon \bfL(\bfp')\lto \bfL(\bfp)$$ which sends
$\vec{l}=\sum_{i=1}^nl_i\vec{x}_i+l\vec{c}$ to
$\phi(\vec{l})=\sum_{i=1}^nl_i\vec{x}_i+l\vec{c}$. Here, $\vec{l}$ is
in its \emph{normal form}, that is, $0\leq l_i<p_i'$ for all $i$. Note
that $\vec{l}\in\bfL(\bfp)$ belongs to the image of $\p$ if and only
if $l_j\ne p_j-1$. Observe that $\phi$ is not a morphism of abelian
groups.

Consider the following functor
$$F\colon \mod^{\bfL(\bfp')}{S(\bfp', \bflambda)} \lto
\mod^{\bfL(\bfp)}{S(\bfp, \bflambda)}$$ which sends
$M=\bigoplus_{\vec{l}\in \bfL(\bfp')} M_{\vec{l}}$ to
$FM=\bigoplus_{\vec{l}\in \bfL(\bfp)} (FM)_{\vec{l}}$, where
$(FM)_{\vec{l}}=M_{\phi^{-1}(\vec{l})}$ if $l_j=0$, and
$(FM)_{\vec{l}}=M_{\phi^{-1}(\vec{l}-\vec{x}_j)}$ otherwise.  The
actions of $u$, $v$ and each $x_i$ on $FM$ are induced from the ones
on $M$, except that $x_j$ acts as the identity
$(FM)_{\vec{l}}\rightarrow (FM)_{ \vec{l}+\vec{x}_j}$ if $l_j=0$.

The functor $F$ identifies the category $\mod^{\bfL(\bfp')}S(\bfp',
\bflambda)$ with the full subcategory of $\mod^{\bfL(\bfp)}S(\bfp,
\bflambda)$ consisting of the modules $N$ such that multiplication
with $x_j$ induces an isomorphism $N_{\vec{l}}\rightarrow
N_{\vec{l}+\vec{x}_j}$ whenever $l_j=0$.

The functor $F$ admits a left adjoint $F_\la$ which sends $N$ to
$F_\la N$ such that $(F_\la N)_{\vec{l}}=N_{\phi(\vec{l})+\vec{x}_j}$.
Similarly, there is a right adjoint $F_\r$ sending $N$ to $F_\r N$
such that $(F_\rho N)_{\vec{l}}=N_{\phi(\vec{l})}$ if $l_j=0$, and
$(F_\r N)_{\vec{l}}=N_{\phi(\vec{l})+\vec{x}_j}$ otherwise. Note that
$x_j$ acts on $(F_\la N)_{\vec{l}}$ as $x_j^2$ if $l_j=p_j-2$, and on
$(F_\rho N)_{\vec{l}}$ as $x_j^2$ if $l_j=0$. We remark that
$F_\r N=F_\la (N(\vec{x}_j))(-\vec{x}_j)$.

All three functors are exact and preserve finite length modules. So
they induce functors between $\coh\bbX'$ and $\coh\bbX$. In
particular, this yields the pursued non-split expansion
$\coh\bbX'\to\coh\bbX$. Note that
$F\mathcal{O}(\vec{l})=\mathcal{O}(\phi(\vec{l}))$ for all $\vec{l}$
in $\bfL(\bfp')$.

Let us remark that the kernel of $F_\la$ on the category of sheaves is
of the form $\add S_\la$ with $S_\la$ a simple sheaf concentrated at
$(x_j)$. Similarly, the kernel of $F_\r$ equals $\add{S_\r}$ with
$S_\r$ another simple sheaf concentrated at $(x_j)$. More precisely,
there is a presentation $0\to \mathcal{O}(-\vec{x}_j)\xto{x_j}
\mathcal{O}\to S_\la\to 0$ and $S_\rho=S_\la(-\vec{x}_j)$. Using the
notation from \ref{se:coh}, we have $S_\la=S_{j1}$ and $S_\r=S_{j, p_j}$.
\end{proof}

\begin{rem}
Let $\bbX = (\bbP^1_k,\bflambda,\bfp)$. The proof shows that the
length $r$ of the filtration of $\A=\coh\bbX$ is determined by the
weight sequence $\bfp$. More precisely, $r=\sum_{i=1}^n(p_i-1)$ and
each category $\A^l$ is of the form $\coh\bbX'$ for some weighted
projective line $\bbX' = (\bbP^1_k,\bflambda,\bfp')$ such that
$p_i'\le p_i$ for $1\le i\le n$.
\end{rem}

\subsection{Vector bundles on weighted projective lines}

Let $\A=\coh\bbX$ be the category of coherent sheaves on a weighted
projective line $\bbX$ and denote by $\vect\bbX$ the full subcategory
formed by all torsion-free objects.  Note that the line bundles in
$\A$ are precisely those objects that are torsion-free and
$1$-simple. The category $\vect\bbX$ admits a Quillen exact structure
such that $\vect \bbX$ has enough projective and injective
objects. Moreover, projective and injective objects coincide; they are
precisely the direct sums of line bundles. Thus $\vect\bbX$ is a
Frobenius category and its stable category $\uvect\bbX$ carries a
triangulated structure; see \cite{KLM} for details.

\begin{thm}
Let $\bbX = (\bbP^1_k,\bflambda,\bfp)$ and $\bbX' =
(\bbP^1_k,\bflambda,\bfp')$ be a pair of weighted projective lines
such that $p_i'\le p_i$ for $1\le i\le n$. Then there is a fully
faithful and exact functor $\coh\bbX'\to\coh\bbX$ that induces the
following recollement
\[ \xymatrixrowsep{3pc} \xymatrixcolsep{2pc}\xymatrix{
\uvect\bbX'\;\ar[rr]|-{F}&&\;\uvect\bbX\; \ar[rr]|-{}
\ar@<1.2ex>[ll]^-{}\ar@<-1.2ex>[ll]_-{}&&
\;\uvect\bbX/\Im F.\ar@<1.2ex>[ll]^-{}\ar@<-1.2ex>[ll]_-{}
}\]
\end{thm}

\begin{proof}
The reduction of weights described in the proof of
Theorem~\ref{th:char} yields a fully faithful functor $i\colon
\coh\bbX'\to\coh\bbX$ which is a composite of $\sum_i(p_i-p_i')$
expansions. This functor has left and right adjoints, and all three
functors preserve vector bundles and line bundles by
Proposition~\ref{pr:vec}. It follows that $i$ restricts to an exact
functor $\vect \bbX'\to\vect\bbX$ which admits exact left and right
adjoints. For the exactness of these functors, one uses Serre duality
and the fact that a sequence $\xi\colon 0\rightarrow
E'\rightarrow E\rightarrow E''\rightarrow 0$ in $\vect\bbX$ is exact
with respect to the specified Quillen exact structure if and only if
$\Hom_\A(L, -)$ sends $\xi$ to an exact sequence for each line bundle
$L$, if and only if $\Hom_\A(-, L)$ sends $\xi$ to an exact sequence
for each line bundle $L$. Thus the induced functors between $\uvect
\bbX'$ and $\uvect\bbX$ form the left hand side of a recollement,
while the right half exists for formal reasons, as explained before in
Lemma~\ref{le:der-adjoint}.
\end{proof}

Not much seems to be known about the right hand term in the
recollement of $\uvect\bbX$. This is in contrast to the recollement of
$\bfD^b(\coh\bbX)$ given in Theorem~\ref{th:char}, where one has a
very explicit description.


\begin{thebibliography}{99}
%
\bibitem{BBD}A. A. Be\u\i linson, J. Bernstein\ and\ P. Deligne,
Faisceaux pervers, in {\it Analysis and topology on singular spaces, I
(Luminy, 1981)}, 5--171, Ast\'erisque, 100, Soc. Math. France, Paris,
1982.
%
\bibitem{BIK} D. Benson, S. B. Iyengar\ and\ H. Krause, Local
  cohomology and support for triangulated categories,
  Ann. Sci. \'Ec. Norm. Sup\'er. (4) {\bf 41} (2008), no.~4, 573--619.
%
\bibitem{CK} X.-W. Chen\ and\ H. Krause, Introduction to coherent sheaves
  on weighted projective lines, arXiv:0911.4473.
%
\bibitem{FGR} R. M. Fossum, P. A. Griffith\ and\ I. Reiten, {\it
  Trivial extensions of abelian categories}, Springer, Berlin, 1975.
%
\bibitem{G}P. Gabriel, Des cat\'egories ab\'eliennes,
Bull. Soc. Math. France {\bf 90} (1962), 323--448.
%
\bibitem{GZ}P. Gabriel\ and\ M. Zisman, {\it Calculus of fractions and
homotopy theory}, Springer-Verlag New York, Inc., New York, 1967.
%
\bibitem{GL1987} W. Geigle\ and\ H. Lenzing, A class of weighted
projective curves arising in representation theory of
finite-dimensional algebras, in {\it Singularities, representation of
algebras, and vector bundles (Lambrecht, 1985)}, 265--297, Lecture
Notes in Math., 1273, Springer, Berlin, 1987.
%
\bibitem{GL1991} W. Geigle\ and\ H. Lenzing, Perpendicular categories
with applications to representations and sheaves, J. Algebra {\bf 144}
(1991), no.~2, 273--343.
%
\bibitem{KLM} D. Kussin, H. Lenzing\ and\ H. Meltzer, Nilpotent
  operators and weighted projective lines, arXiv:1002.3797.
%
\bibitem{L} H. Lenzing, Hereditary Noetherian categories with a
tilting complex, Proc. Amer. Math. Soc. {\bf 125} (1997), no.~7,
1893--1901.
%
\bibitem{Le1998} H. Lenzing, Representations of finite-dimensional
  algebras and singularity theory, in {\it Trends in ring theory
    (Miskolc, 1996)}, 71--97, Amer. Math. Soc., Providence, RI, 1998.
%
\bibitem{RvB} I. Reiten\ and\ M. Van den Bergh, Noetherian hereditary
  abelian categories satisfying Serre duality,
  J. Amer. Math. Soc. {\bf 15} (2002), no.~2, 295--366.
%
\bibitem{Ri1980} C. M. Ringel, On algorithms for solving vector space
  problems. II. Tame algebras, in {\it Representation theory, I
    (Proc. Workshop, Carleton Univ., Ottawa, Ont., 1979)}, 137--287,
  Lecture Notes in Math., 831, Springer, Berlin, 1980.
%
\bibitem{Sch1985} A. H. Schofield, {\it Representation of rings over
  skew fields}, Cambridge Univ. Press, Cambridge, 1985.
%
\bibitem{Sch} A. Schofield, Generic representations of quivers,
  unpublished manuscript.
%
\end{thebibliography}
\end{document}